\documentclass[12pt]{article}
\usepackage{amsmath}
\usepackage{amsthm}
\usepackage{amsfonts}
\usepackage{amssymb}
\usepackage{graphicx}
\usepackage{xypic}
\usepackage{a4wide}
\usepackage{enumerate}
\usepackage{mathrsfs}
\newtheorem{theorem}{Theorem}[section]
\newtheorem{proposition}[theorem]{Proposition}
\newtheorem{corollary}[theorem]{Corollary}
\newtheorem{lemma}[theorem]{Lemma}
\theoremstyle{definition}

\newtheorem{definition}[theorem]{Definition}

\newcommand{\bcomod}[2]{{}^{#1}\mathcal{M}^{#2}}
\newcommand{\cat}[1]{\mathbf{#1}}
\newcommand{\cohom}[3]{\mathrm{h}_{#1}(#2,#3)}
\newcommand{\coring}[1]{\mathfrak{#1}}
\newcommand{\cotensor}[1]{\square_{#1}}
\renewcommand{\hom}[3]{\mathrm{Hom}_{#1}(#2,#3)}
\newcommand{\lcomod}[1]{{}^{#1}\mathcal{M}}
\newcommand{\rcomod}[1]{\mathcal{M}^{#1}}
\newcommand{\rmod}[1]{\mathcal{M}_{#1}}
\newcommand{\lmod}[1]{_{#1}\mathcal{M}}

\newcommand{\tensor}[1]{\otimes_{#1}}
\newcommand{\e}[2]{\mathrm{e}_{#1}(#2)}

\newcommand{\pic}[1]{\mathrm{Pic}_k(#1)}
\newcommand{\rpic}[1]{\mathrm{Pic}_k^r{(#1)}}
\newcommand{\lpic}[1]{\mathrm{Pic}_k^l{(#1)}}
\newcommand{\aut}[1]{\mathrm{Aut}_k(#1)}
\newcommand{\inn}[1]{\mathrm{Inn}_k(#1)}
\newcommand{\rinn}[1]{\mathrm{Inn}_k^r{(#1)}}
\newcommand{\out}[1]{\mathrm{Out}_k(#1)}
\newcommand{\rout}[1]{\mathrm{Out}_k^r{(#1)}}
\begin{document}

\title{The Picard Group of Corings}
\author{Mohssin Zarouali-Darkaoui}
\date{\empty}

\maketitle

\section*{Introduction}

Corings were introduced by Sweedler in \cite{Sweedler:1975}. A
coring over an associative algebra with unit over a commutative ring
$k$, $A$, is an $A$-bimodule $\coring{C}$ with two $A$-bimodule maps
$\Delta_{\coring{C}}:\coring{C}\rightarrow\coring{C}\tensor{A}{\coring{C}}$
(comultiplication) and $\epsilon_{\coring{C}}:\coring{C}\rightarrow
A$ (counit) such that the same diagrams as for coalgebras are
commutative. Recently, corings were intensively studied. The main
motivation of this studies is the observation of Takeuchi, Theorem
4.1. For a detailed study of corings, we refer to
\cite{Brzezinski/Wisbauer:2003}.

The connection between the Picard group of Azumaya algebras and its
automorphisms were studied in \cite{Rosenberg/Zelinsky:1961}. Bass
generalized these connections for arbitrary algebras, see
\cite{Bass:1968}. Further study of the Picard group of algebras is
given in \cite{Frohlich:1973}. In \cite{Beattie/Delrio:1996}, the
authors have given the versions of these connections for rings with
local units, in order to study the Picard group of the category
$R-gr$, where $R$ is a $G$-graded ring. In
\cite{Torrecillas/Zhang:1996}, the authors studied these connections
for coalgebras over fields.

The purpose of this note is to introduce and study the right Picard
group of corings. The motivation is the fact that there is an
isomorphism of groups between the Picard group of the category
$\rcomod{\coring{C}}$, for a certain coring $\coring{C}$, and the
right Picard group of $\coring{C}$ which is defined as the group of
the isomorphism classes of right invertible $\coring{C}$-bicomodules
(=$T(\coring{C})$ for some $k$-autoequivalence of
$\rcomod{\coring{C}}$, $T$) with the composition law induced by the
cotensor product (see Proposition \ref{picard-cat-crg}). We extend
the exact sequence $1\to \inn{-}\to \aut{-}\to \pic{-}$ from
algebras and coalgebras over fields to corings. We also extend a
result which is useful to show that a given coring $\coring{C}$ have
the right Aut-Pic property, i.e. where the morphism
$\aut{\coring{C}}\to \rpic{\coring{C}}$ is an epimorphism (see
Proposition 2.7). Of course, our Aut-Pic property for corings
extends that of algebras \cite{Bolla:1984}, and that of coalgebras
over fields \cite{Cuadra/Garcia/Torrecillas:2000}. All of the
examples of algebras and coalgebras having the Aut-Pic property that
are given in \cite{Bolla:1984,Cuadra/Garcia/Torrecillas:2000}, are
corings having this property. In this note we give some new examples
of corings having the Aut-Pic property. We also simplify the
computation of the right Picard group of several interesting corings
(see Proposition \ref{examples}). Finally, in section 4, we give the
corresponding exact sequences for the category of entwined modules
over an entwining structure, the category of Doi-Koppinen-Hopf
modules over a Doi-Koppinen structure, and the category of graded
modules by a $G$-set, where $G$ is a group.

\section{Preliminaries}

Throughout this note and unless otherwise stated, $k$ denote a
commutative ring (with unit), $A,$ $A',$ $A'',$ $A_1,$ $A_2,$ $B,$
$B_1,$ and $B_2$ denote associative and unitary algebras over $k$,
and $\coring{C},$ $\coring{C}',$ $\coring{C}'',$ $\coring{C}_1,$
$\coring{C}_2,$ $\coring{D},$ $\coring{D}_1$ and $\coring{D}_2$
denote corings over $A,$ $A',$ $A'',$ $A_1,$ $A_2,$ $B,$ $B_1,$ and
$B_2$, respectively. The notation $\otimes$ will stand for the
tensor product over $k$.

A category $\cat{C}$ is said to be $k$-\emph{category} if for every
$M$ and $N$ in $\cat{C}$, $\hom{\cat{C}}{M}{N}$ is a $k$-module, and
the composition is $k$-bilinear. An abelian category which is a
$k$-category is said to be $k$-\emph{abelian category}. A functor
between $k$-categories is said to be $k$-\emph{functor} or
$k$-\emph{linear functor} if it is $k$-linear on the $k$-modules of
morphisms. A functor between $k$-categories is said to be a
$k$-\emph{equivalence} if it is $k$-linear and an equivalence.

We recall from \cite{Sweedler:1975} that an $A$-\emph{coring} is an
$A$-bimodule $\coring{C}$ with two $A$-bimodule maps $\Delta :
\coring{C} \rightarrow \coring{C} \tensor{A} \coring{C}$ and
$\epsilon : \coring{C} \rightarrow A$ such that $(\coring{C}
\tensor{A} \Delta) \circ \Delta = (\Delta \tensor{A} \coring{C})
\circ \Delta$ and $(\epsilon \tensor{A} \coring{C}) \circ \Delta = (
\coring{C} \tensor{A} \epsilon) \circ \Delta = 1_\coring{C}$.
$\coring{C}=A$ endowed with the obvious structure maps is an
$A$-coring. A \emph{right} $\coring{C}$-\emph{comodule} is a pair
$(M,\rho_M)$ where $M$ is a right $A$-module and $\rho_M: M
\rightarrow M \tensor{A} \coring{C}$ (coaction) is an $A$-linear map
satisfying $(M \tensor{A} \Delta) \circ \rho_M = (\rho_M \tensor{A}
\coring{C}) \circ \rho_M$, and $(M \tensor{A} \epsilon) \circ \rho_M
= 1_M$. A \emph{morphism} of right $\coring{C}$-comodules
$(M,\rho_M)$ and $(N,\rho_N)$ is a right $A$-linear map $f: M
\rightarrow N$ such that $(f \tensor{A} \coring{C}) \circ \rho_M =
\rho_N \circ f$. Right $\coring{C}$-comodules with their morphisms
form a $k$-category $\rcomod{\coring{C}}$. Coproducts and cokernels
(and then inductive limits) in $\rcomod{\coring{C}}$ exist and they
coincide respectively with coproducts and cokernels in the category
of right $A$-modules $\rmod{A}$. If ${}_A\coring{C}$ is flat, then
$\rcomod{\coring{C}}$ is a $k$-abelian category. Moreover it is a
Grothendieck category. When $\coring{C}=A$ is the obvious
$A$-coring, $\rcomod{A}$ is the category of right $A$-modules
$\rmod{A}$.

Now assume that the $A'-A$-bimodule $M$ is also a left comodule over
an $A'$-coring $\coring{C}'$ with structure map $\lambda_M : M
\rightarrow \coring{C}' \tensor{A'} M$. Assume moreover that
$\rho_M$ is $A'$-linear, and $\lambda_M$ is $A$-linear. It is clear
that $\rho_M : M \rightarrow M \tensor{A} \coring{C}$ is a morphism
of left $\coring{C}'$-comodules if and only if $\lambda_M : M
\rightarrow \coring{C}' \tensor{A'} M$ is a morphism of right
$\coring{C}$-comodules. In this case, we say that $M$ is a
$\coring{C}'-\coring{C}$-bicomodule. A morphism of bicomodules is a
morphism of right and left comodules. Then we obtain a $k$-category
$^\coring{C'}\mathcal{M}^\coring{C}$. If $_A\coring{C}$ and
$\coring{C'}_{A'}$ it is a Grothendieck category. If
$\coring{C}'=A', \coring{C}=A$, then
$^\coring{C'}\mathcal{M}^\coring{C}$ is the category of
$A'-A$-bimodules, $_{A'}\mathcal{M}_A$.

A coring $\coring{C}$ is said to be \emph{coseparable} if the
comultiplication map $\Delta_{\coring{C}}$ is a section in the
category $^\coring{C}\mathcal{M}^\coring{C}$. Obviously the trivial
$A$-coring $\coring{C} = A$ is coseparable.

Let $Z$ be a left $A$-module and $f: X \rightarrow Y$ a morphism in
$\rmod{A}$. Following \cite[40.13]{Brzezinski/Wisbauer:2003} we say
that $f$ is \emph{$Z$-pure} when the functor $-\tensor{A} Z$
preserves the kernel of $f$. If $f$ is $Z$-pure for every $Z \in
{}\lmod{A}$ then we say simply that $f$ is \emph{pure} in
$\rmod{A}$.

Let $f:M\to N$ be a morphism in $\bcomod{\coring{C}'}{\coring{C}}$,
and let $\operatorname{ker}(f)$ be its kernel in
$_{A'}\mathcal{M}_A$. It is easy to show that if
$\operatorname{ker}(f)$ is $\coring{C}'_{A'}$-pure and
$_A\coring{C}$-pure, and the following
\begin{equation*}
\operatorname{ker}(f)\otimes_A\coring{C} \otimes_A\coring{C}, \quad
\coring{C}'\otimes_{A'}\coring{C}'\otimes_{A'}\operatorname{ker}(f)
\quad \textrm{and} \quad \coring{C}'\otimes_{A'}\operatorname{ker}
(f)\otimes_A\coring{C}
\end{equation*}
are injective maps, then $\operatorname{ker}(f)$ is the kernel of
$f$ in ${}^{\coring{C}'}\mathcal{M}^{\coring{C}}$. This is the case
if $f$ is $(\coring{C}'\otimes_{A'}\coring{C}')_{A'}$-pure,
$_A(\coring{C} \otimes_A\coring{C})$-pure, and
$\coring{C}'\otimes_{A'}f$ is $_A\coring{C}$-pure (e.g. if
$\coring{C}'_{A'}$ and $_A\coring{C}$ are flat, or if $\coring{C}$
is a coseparable $A$-coring).

Now let
 $M\in{}\bcomod{\coring{C}'}{\coring{C}}$ and
$N\in{}^{\coring{C}}\mathcal{M}^{\coring{C}''}$. The map
\[
\omega_{M,N}=\rho_M\otimes_AN-M\otimes_A\lambda_N:M\otimes_AN
\rightarrow M\otimes_A\coring{C}\otimes_AN
\]
is a $\coring{C}'-\coring{C}''$-bicomodule map. Its kernel in
$_{A'}\mathcal{M}_{A''}$ is the \emph{cotensor product} of $M$ and
$N$, and it is denoted by $M\square_{\coring{C}}N$ . From the above
consideration, if $\omega_{M,N}$ is
$(\coring{C}'\otimes_{A'}\coring{C}')_{A'}$-pure,
$_{A''}(\coring{C}'' \otimes_{A''}\coring{C}'')$-pure, and
$\coring{C}'\otimes_{A'}\omega_{M,N}$ is $_{A''}\coring{C}''$-pure
(e.g. if $\coring{C}'_{A'}$ and $_{A''}\coring{C}''$ are flat, or if
$\coring{C}$ is a coseparable $A$-coring), then
$M\square_{\coring{C}}N$ is the kernel of $\omega_{M,N}$ in
${}^{\coring{C}'}\mathcal{M}^{\coring{C}''}$.

If for every $M\in{}\bcomod{\coring{C}'}{\coring{C}}$ and
$N\in{}^{\coring{C}}\mathcal{M}^{\coring{C}''}$, $\omega_{M,N}$ is
$\coring{C}'_{A'}$-pure and $_{A''}\coring{C}''$-pure, then we have
a $k$-linear bifunctor
\begin{equation}\label{cotensorbifunctor}
\xymatrix{-\cotensor{\coring{C}}-:{}^{\coring{C}'}\mathcal{M}^{\coring{C}}
\times{}^{\coring{C}}\mathcal{M}^{\coring{C}''}\ar[r]&
{}^{\coring{C}'}\mathcal{M}^{\coring{C}''}}.
\end{equation} If in particular $\coring{C}'_{A'}$ and
$_{A''}\coring{C}'' $ are flat, or if $\coring{C}$ is a coseparable
$A$-coring, then the bifunctor \eqref{cotensorbifunctor} is well
defined. In the special case when $\coring{C}=A$,
$-\cotensor{\coring{C}}-=-\tensor{A}-$.

Throughout this note, for all coring $\coring{C}$ over $A$, we have
$_A\coring{C}$ and $\coring{C}_A$ are flat.

\section{The Picard group of corings}

A \emph{coring homomorphism} \cite{Gomez:2002} from the coring
$(\coring{C}:A)$ to the coring $(\coring{D}:B)$ is a pair
$f=(\varphi,\rho)$, where $\rho:A\rightarrow B$ is a homomorphism of
$k$-algebras and $\varphi:\coring{C}\rightarrow\coring{D}$ is a
homomorphism of $A$-bimodules such that
\[
\epsilon_{\coring{D}}\circ\varphi=\rho\circ\epsilon_{\coring{C}}\qquad\textrm{and
}\qquad\Delta_{\coring{D}}\circ\varphi=\omega_{\coring{D},\coring{D}}\circ(
\varphi\otimes_A\varphi)\circ\Delta_{\coring{C}},
\]
where $\omega_{\coring{D},\coring{D}}:\coring{D}\otimes_A\coring{D}%
\rightarrow\coring{D}\otimes_B\coring{D}$ is the canonical map
induced by $\rho:A\rightarrow B.$ Corings over $k$-algebras and
theirs morphisms form a category. We denote it by $\mathbf{Crg}_k$.
Notice that the category of $k$-algebras, $\mathbf{Alg}_k$, and the
category of $k$-coalgebras, $\mathbf{Coalg}_k$, are full
subcategories of $\mathbf{Crg}_k$. Notice also that
$(\varphi,\rho)\in \mathbf{Crg}_k$ is an isomorphism if and only if
$\varphi$ and $\rho$ are bijective.

\medskip
To state the definition of the right Picard group of a coring we
need to recall the following result.

\begin{proposition}\label{equivalence}(\cite[Proposition
3.1]{Zarouali:2005}) Suppose that $_A\coring{C}$, $\coring{C}_A$,
$_B\coring{D}$ and $\coring{D}_B$ are flat. Let
$X\in{}^{\coring{C}}\mathcal{M}^{\coring{D}}$ and
$\Lambda\in{}^{\coring{D}}\mathcal{M}^{\coring{C}}$. \item The
following statements are equivalent:
\begin{enumerate}[(1)]
\item $(-\square _{\coring{C}}X,-\square_{\coring{D}}\Lambda)$ is a
pair of inverse equivalences;
\item there exist bicomodule isomorphisms
\[
 f:X\square_{\coring{D}}\Lambda \rightarrow\coring{C}\quad
\text{and}\quad g:\Lambda\square_{\coring{C}}X\rightarrow\coring{D}
\]
in $^\coring{C}\mathcal{M}^\coring{C}$ and $^{\coring{D}}%
\mathcal{M}^{\coring{D}}$ respectively, such that
\begin{enumerate}[(a)] \item $_A X$ and $_B \Lambda$ are flat,
and
$\omega_{X,\Lambda}=\rho{}_X\otimes_B\Lambda-X\otimes_A\rho{}_\Lambda$
is pure in $\lmod{A}$ and
$\omega_{\Lambda,X}=\rho{}_\Lambda\otimes_AX-\Lambda\otimes_B\rho{}_X$
is pure in $\lmod{B}$, or \item ${}_{\coring{C}}X$ and
${}_{\coring{D}}\Lambda$ are coflat.
\end{enumerate}
\end{enumerate}
In such a case the diagrams
\begin{equation*}\label{unitcounit}
\xymatrix{\Lambda\cotensor{\coring{C}}X\cotensor{\coring{D}}\Lambda
\ar[rr]^{\Lambda\cotensor{\coring{C}}f}
\ar[d]_{g\cotensor{\coring{D}}\Lambda} &&
\Lambda\cotensor{\coring{C}}\coring{C} \ar[d]^\simeq
\\\coring{D}\cotensor{\coring{D}}\Lambda \ar[rr]^\simeq &&
\Lambda}\qquad
\xymatrix{X\cotensor{\coring{D}}\Lambda\cotensor{\coring{C}}X\ar[rr]^{f\cotensor{\coring{C}}X}
\ar[d]_{X\cotensor{\coring{D}}g}&&
\coring{C}\cotensor{\coring{C}}X\ar[d]^\simeq
\\X\cotensor{\coring{D}}\coring{D} \ar[rr]^\simeq && X}
\end{equation*}
commute.

If $A$ and $B$ are von Neumann regular rings, or if $\coring{C}$ and
$\coring{D}$ are coseparable corings (without the condition
``$_A\coring{C}$, $\coring{C}_A$, $_B\coring{D}$ and $\coring{D}_B$
are flat''), the conditions (a) and (b) can be deleted.
\end{proposition}

\medskip

Let $\cat{C}$ be a $k$-abelian category, $\pic{\cat{C}}$
\cite{Bass:1968} is the group of isomorphism classes $(T)$ of
$k$-equivalences $T:\cat{C}\to\cat{C}$ with the composition law
$(T)(S)=(ST)$.

\begin{definition}
We say that a $\coring{C}-\coring{D}$-bicomodule $X$ is \emph{right
invertible} if the functor
$-\cotensor{\coring{C}}X:\rcomod{\coring{C}}\to\rcomod{\coring{D}}$
is an equivalence. We say that $\coring{C}$ and $\coring{D}$ are
\emph{Morita-Takeuchi right equivalent} if there is a right
invertible $\coring{C}-\coring{D}$-bicomodule.
\end{definition}

The isomorphism classes $(X)$ of right invertible
$\coring{C}$-bicomodules with the composition law
$$(X_1)(X_2)=(X_1\cotensor{\coring{C}}X_2),$$
form the \emph{right Picard group} of $\coring{C}$. We denote it by
$\rpic{\coring{C}}$. This law is well defined since
$-\cotensor{\coring{C}}(X_1\cotensor{\coring{C}}X_2)\simeq
(-\cotensor{\coring{C}}X_1)\cotensor{\coring{C}}X_2$ (From
\cite[22.5, 22.6(iii)]{Brzezinski/Wisbauer:2003} and
$_\coring{C}X_2$ is coflat). The associativity of this law follows
from the associativity of the cotensor product ($_\coring{C}X_1$ and
$_\coring{C}X_2$ are coflat). The isomorphism class $(\coring{C})$
is the identity element. $(X)^{-1}=(\Lambda)$, where $\Lambda$ is
such that $(-\cotensor{\coring{C}}X,-\cotensor{\coring{C}}\Lambda)$
is a pair of inverse equivalences.

It follows from \cite[Theorems 3.8, 3.13]{Zarouali:2005} that if
$\coring{C}$ is coseparable (this case includes that of algebras),
or cosemisimple, or $A$ is von Neumann regular ring, or if
$\coring{C}$ is a coalgebra over a QF ring, then a
$\coring{C}$-bicomodule $X$ is right invertible if and only if it is
left invertible. Hence, for each case we have
$\rpic{\coring{C}}=\lpic{\coring{C}}$. We will denote this last by
$\pic{\coring{C}}$.

\medskip

The following proposition is the motivation of the study of the
Picard group of corings. The proof is straightforward using
\cite[Theorem 2.3]{Zarouali:2004}.

\begin{proposition}\label{picard-cat-crg}
There are inverse isomorphisms of groups
$$\xymatrix{\rpic{\coring{C}}\ar@<1ex>[r]^\alpha & \pic{\rcomod{\coring{C}}}
\ar@<1ex>[l]^\beta},$$ $\alpha((X))=(-\cotensor{\coring{C}}X)$ and
$\beta((T))=(T(\coring{C}))$.
\end{proposition}

\medskip
Now we introduce two categories, $\mathfrak{Crg}_k$ and
$\mathfrak{MT-Crg}_k^r$. The objects of both are the corings. The
morphisms of $\mathfrak{Crg}_k$ are the isomorphisms of corings. The
morphism of $\mathfrak{MT-Crg}_k^r$ are the Morita-Takeuchi right
equivalences. A \emph{right Morita-Takeuchi equivalence}
$\coring{C}\sim\coring{D}$ is simply an isomorphism class $(M)$
where $_\coring{D}M_\coring{C}$ is a right invertible
$\coring{D}-\coring{C}$-bicomodule. If $(N)$ is a Morita-Takeuchi
right equivalence $\coring{D}\sim\coring{C}'$, the composite is
given by $(N\cotensor{\coring{D}}M)$.

The categories $\mathfrak{Crg}_k$ and $\mathfrak{MT-Crg}_k^r$ are
groupoids (i.e., categories in which every morphism is an
isomorphism). Then $\aut{\coring{C}}$ and $\rpic{\coring{C}}$ are
respectively the endomorphism group of the object $\coring{C}$ in
$\mathfrak{Crg}_k$ and in $\mathfrak{MT-Crg}_k$.

Let $f:(\coring{C}:A)\to(\coring{D}:B)$ and
$g:(\coring{C}_1:A_1)\to(\coring{D}_1:B_1)$ be two morphisms of
corings, and let $_{\coring{C}_1}X_{\coring{C}}$ be a bicomodule.

We consider the right induction functor $-\tensor{A}B:
\rcomod{\coring{C}}\rightarrow \rcomod{\coring{D}}$ defined in
\cite{Gomez:2002}. The right coaction on the right $B$-module
$M\tensor{A}B$ is given by:
$$\rho_{M\tensor{A}B}:M\tensor{A}B\rightarrow
M\tensor{A}B\tensor{B}\coring{D}\simeq M\tensor{A}\coring{D},\;
m\tensor{A}b\mapsto \sum m_{(0)}\tensor{A}\varphi(m_{(1)})b,$$ where
$\rho_M(m)=\sum m_{(0)}\tensor{A}m_{(1)}$. By a similar way we
define the left induction functor $B_1\tensor{A_1}-:
\lcomod{\coring{C_1}}\rightarrow \lcomod{\coring{D_1}}$.

Now consider $_gX_f:=B_1\tensor{A_1}X\tensor{A}B$ which is a
$\coring{D}_1-\coring{D}$-bicomodule, and $X_f:=X\tensor{A}B$ which
is a $\coring{C}_1-\coring{D}$-bicomodule. It can be showed easily
that $_1X_f\simeq X_f$ as bicomodules, and that if
$\xymatrix{\coring{C}\ar[rr]^{f=(\varphi,\rho)} && \coring{D}
\ar[rr]^{f'=(\varphi',\rho')}&& \coring{C}'}$ are morphisms of
corings, then $(X_f)_{f'}\simeq X_{f'f}$ as bicomodules. Moreover,
from \cite[Theorem 2.4]{Zarouali:2004}, $X_f\simeq
X\cotensor{\coring{C}}\coring{C}_f$ as
$\coring{C}_1-\coring{D}$-bicomodules. Finally, it is easy to show
that if $f$ is an isomorphism, then $X_f\simeq (X_f)'$ as
bicomodules, where $(X_f)'=X_f$, and the left
$\coring{C}_1$-comodule structure on $(X_f)'$ is the same as for
$X_f$. The right $B$-module structure on $(X_f)'$ is given by
$xb=x\rho^{-1}(b)$, for $x\in X,b\in B$. The right
$\coring{D}$-comodule structure on it is given by
$\rho_{(X_f)'}(x)=\sum x_{(0)}\tensor{B}\varphi(x_{(1)})$, where
$x\in X$ and $\rho_X(x)=\sum x_{(0)}\tensor{A}x_{(1)}$.

\begin{lemma}\label{functor}
We have a functor $$\Omega:\mathfrak{Crg}_k \to
\mathfrak{MT-Crg}_k^r,$$ which is the identity on objects and
associates with an isomorphism of corings
$f:(\coring{C}:A)\to(\coring{D}:B)$ the isomorphism class of the
right invertible $\coring{D}-\coring{C}$-bicomodule $_f\coring{C}$.
\end{lemma}

\begin{proof}
Let $f=(\varphi,\rho):(\coring{C}:A)\to(\coring{D}:B)$ be an
isomorphism of corings. Since
$B\tensor{A}\coring{C}\tensor{A}B\to\coring{D},\;
b\tensor{A}c\tensor{A}b\mapsto b\varphi(c)b$ is a morphism of
$B$-corings (see \cite[24.1]{Brzezinski/Wisbauer:2003}), it follows
from \cite[Theorem 4.1]{Zarouali:2005} that $_f\coring{C}$ is a
right invertible $\coring{D}-\coring{C}$-bicomodule. Now let
$\xymatrix{\coring{C}\ar[r]^f & \coring{D} \ar[r]^g& \coring{C}'}$
be two morphisms in $\mathfrak{Crg}_k$. Then
$_{gf}\coring{C}\simeq{}_g(_f\coring{C})\simeq {}_g\coring{D}
\cotensor{\coring{D}}{}_f\coring{C}$ is an isomorphism of
$\coring{C'}-\coring{C}$-bicomodules.
\end{proof}

Now we are ready to state and prove the main result of this section.
\begin{theorem}\label{coringexseq}
Let $\coring{C}$ be an $A$-coring such that $_A\coring{C}$ and
$\coring{C}_A$ are flat. \begin{enumerate}[(1)]
\item We have an exact sequence
$$\xymatrix{1\ar[r] & \rinn{\coring{C}}\ar[r]&
\aut{\coring{C}}\ar[r]^\Omega & \rpic{\coring{C}}},$$ where
$\rinn{\coring{C}}$ is the set of $f=(\varphi,\rho)\in
\aut{\coring{C}}$ such that there is $p\in \coring{C}^*$ invertible
with $$\sum \varphi(c_{(1)})p(c_{(2)})=\sum p(c_{(1)})c_{(2)},$$ for
every $c\in \coring{C}$. As for algebras and coalgebras over fields,
we call $\rinn{\coring{C}}$ the \emph{right group of inner
automorphisms} of $\coring{C}$

In particular there is a monomorphism from the quotient group
$\rout{\coring{C}}:=\aut{\coring{C}}/\rinn{\coring{C}}$ to
$\rpic{\coring{C}}$. As for algebras and coalgebras over fields, we
call $\rout{\coring{C}}$ the \emph{right outer automorphism group}
of $\coring{C}$.
\item If $\coring{C}$ and $\coring{D}$ are Morita-Takeuchi right
equivalent, then $\rpic{\coring{C}}\simeq\rpic{\coring{D}}$.
\end{enumerate}
\end{theorem}

\begin{proof}
$(1)$ Let $f=(\varphi,\rho)\in \aut{\coring{C}}$ and
$_f\coring{C}\simeq \coring{C}$ as $\coring{C}$-bicomodules.
 We have $_f\coring{C}\simeq(_f\coring{C})'$ as
$\coring{C}$-bicomodules (see the remark just before Lemma
\ref{functor}). Now let $h:(_f\coring{C})'\to \coring{C}$ be an
 isomorphism of $\coring{C}$-bicomodules.
$h$ is a morphism of right $\coring{C}$-bicomodules if and only if
there is $p\in \coring{C}^*$ invertible such that $h(c)=\sum
p(c_{(1)})c_{(2)},$ for every $c\in \coring{C}$ (see
\cite[18.12(1)]{Brzezinski/Wisbauer:2003}). In such a case,
$p=\epsilon h$. $h$ is a morphism of left $\coring{C}$-bicomodules
means that $$h(ac)=\rho(a)h(c),\quad(*)$$ for every $a\in A, c\in
\coring{C}$, and for every $c\in \coring{C}$,
$$\sum\varphi(c_{(1)})\tensor{A}h(c_{(2)})=\sum
h(c)_{(1)}\tensor{A}h(c)_{(2)}.\quad(**)$$ The condition $(**)$ is
equivalent to
$$\sum \varphi(c_{(1)})\tensor{A}p(c_{(2)})c_{(3)}=\sum
p(c_{(1)})c_{(2)}\tensor{A}c_{(3)},$$ for every $c\in \coring{C}$.
Hence
$$\sum \varphi(c_{(1)})p(c_{(2)})=\sum p(c_{(1)})c_{(2)},$$ for every
$c\in \coring{C}$. Since $\varphi$ is left $A$-linear, $(*)$ holds.
The converse is obvious.

$(2)$ It follows immediately from that $\coring{C}$ and $\coring{D}$
are isomorphic to each other in the category
$\mathfrak{MT-Crg}_k^r$. More explicitly, Let
$_\coring{C}M_\coring{D}$ be a right invertible bicomodule. The map
$$(X)\mapsto (M^{-1}\cotensor{\coring{C}}X\cotensor{\coring{C}}M)$$
is an isomorphism from $\rpic{\coring{C}}$ to $\rpic{\coring{D}}$.
\end{proof}

The last theorem gives a well-know result of Bass (see
\cite[Proposition II (5.2)(3)]{Bass:1968} or \cite[Theorems 55.9,
55.11]{Curtis/Reiner:1987}), and a generalization to the case of
coalgebras over rings of the particular case of \cite[Theorem
2.7]{Torrecillas/Zhang:1996} where $R=k$.

\begin{corollary}
\begin{enumerate}[(1)]
\item For a $k$-algebra $A$, we have an exact sequence
$$\xymatrix{1\ar[r] & \inn{A}\ar[r]& \aut{A}\ar[r]^\Omega & \pic{A}},$$
where $\inn{A}:=\{a\mapsto bab^{-1}\mid b \;\textrm{invertible in}
\;A\},$ the \emph{group of inner automorphisms} of $A$.
\item  For a $k$-coalgebra $C$ such that $_kC$ is flat, we have an
exact sequence
$$\xymatrix{1\ar[r] & \inn{C}\ar[r]& \aut{C}\ar[r]^\Omega &
\rpic{C}},$$ where $\inn{C}$ (the \emph{group of inner
automorphisms} of $C$) is the set of $\varphi\in \aut{C}$ such that
there is $p\in C^*$ invertible with $\varphi(c)=\sum
p(c_{(1)})c_{(2)}p^{-1}(c_{(3)}),$ for every $c\in C.$

\end{enumerate}
\end{corollary}

\begin{proof}
$(1)$ It suffices to see that for $p\in \coring{C}^*$ where
$\coring{C}=A$, $p$ is invertible if and only if
$p=\lambda_b:a\mapsto ba$ and $b$ is invertible in $A$.

$(2)$ Suppose that $\varphi\in \aut{C}$ and $\sum
\varphi(c_{(1)})p(c_{(2)})=\sum p(c_{(1)})c_{(2)},$ for every $c\in
C$. Then $\varphi(c)=\sum\varphi(c_{(1)})\epsilon_C(c_{(2)})=\sum
\varphi(c_{(1)})p(c_{(2)})p^{-1}(c_{(3)})= \sum
p(c_{(1)})c_{(2)}p^{-1}(c_{(3)}),$ for every $c\in C$. The converse
is obvious.
\end{proof}

\medskip
In order to consider the one-sided comodule structure of invertible
bicomodules we need

\begin{proposition}\label{oneside}
Let $(M),(N)\in \rpic{\coring{C}}$. Then
\begin{enumerate}[(1)] \item $(N)\in \mathrm{Im}\Omega.(M)$  if and
only if $N\simeq {}_fM$ as bicomodules for some $f\in
\aut{\coring{C}}$. \item $_AM_\coring{C}\simeq {}_AN_\coring{C}$ if
and only if $N\simeq {}_fM$ as bicomodules for some
$f=(\varphi,1_A)\in \aut{\coring{C}}$.
\end{enumerate}
\end{proposition}

\begin{proof}
$(1)$ Straightforward from the fact that for every $f\in
\aut{\coring{C}}$, there is an isomorphism of
$\coring{C}$-bicomodules
$_fM\simeq{}_f\coring{C}\cotensor{\coring{C}}M$.

$(2)$ $(\Leftarrow)$ Obvious.  $(\Rightarrow)$ let
$h:{}_AM_\coring{C}\to {}_AN_\coring{C}$ be a bicomodule
isomorphism. Since $M$ is quasi-finite as a right
$\coring{C}$-comodule (see \cite[Proposition 3.6]{Zarouali:2005}),
it has a structure of $\e{\coring{C}}{M}-\coring{C}$-bicomodule.
Then $N$ has a structure of
$\e{\coring{C}}{M}-\coring{C}$-bicomodule induced by $h$, and $h$ is
an $\e{\coring{C}}{M}-\coring{C}$-bicomodule isomorphism. On the
other hand, $N$ is quasi-finite as a right $\coring{C}$-comodule.
Let $F=:\cohom{\coring{C}}{N}{-}:\rcomod{\coring{C}}\to
\rcomod{\e{\coring{C}}{M}}$ be the cohom functor, and let
$\theta:1_{\rcomod{\coring{C}}}\to
F(-)\cotensor{\e{\coring{C}}{M}}N$ and
$\chi:F(-\cotensor{\e{\coring{C}}{M}}N)\to
1_{\rcomod{\e{\coring{C}}{M}}}$ be respectively the unit and the
counit of this adjunction. From \cite[Proposition
5.2]{ElKaoutit/Gomez:2003} and its proof, the map
$\varphi:=\chi_{\e{\coring{C}}{M}}\circ
F(\lambda_N^h):\e{\coring{C}}{N}\to \e{\coring{C}}{M}$, where
$\lambda_N^h:N\to \e{\coring{C}}{M}\cotensor{\e{\coring{C}}{M}}N$ is
the left comodule structure map on $N$, is a morphism of
$A$-corings, and making commutative the diagram
$$\xymatrix{N
\ar[ddrr]^{\lambda_N^h}\ar[rr]^{\theta_N} &&
\e{\coring{C}}{N}\tensor{A}N\ar[dd]^{\varphi\tensor{A}1_N}
\\ &&&\\
& & \e{\coring{C}}{M}\tensor{A}N.}$$  Since $F$ is a left adjoint to
$-\cotensor{\e{\coring{C}}{M}}N$ which is an equivalence,
$\chi_{\e{\coring{C}}{M}}$ and $\varphi$ are isomorphisms. Since $M$
and $N$ are right invertibles,
$\e{\coring{C}}{M}\simeq\coring{C}\simeq\e{\coring{C}}{N}$ (see
\cite[Proposition 3.6]{Zarouali:2005}). We identify
$\e{\coring{C}}{M}$ and $\e{\coring{C}}{N}$ with $\coring{C}$. Then
$\varphi$ is a coring endomorphism of $\coring{C}$. Hence $M\simeq
{}_fN$ as $\coring{C}$-bicomodules, where $f=(\varphi,1_A)$.
\end{proof}

As a immediate consequence of the last result, we obtain a
generalization to the case of coalgebras over rings of the
particular case of \cite[Proposition 2.8]{Torrecillas/Zhang:1996}
where $R=k$. Its version for algebras is due to Bass, see
\cite[Proposition II (5.2)(4)]{Bass:1968} or \cite[Theorem
55.12]{Curtis/Reiner:1987}.

\begin{corollary}
Let $C$ be a $k$-coalgebra such that $_kC$ is flat, and let
$(M),(N)\in \rpic{C}$. Then $M_C\simeq N_C$ if and only if $(N)\in
\mathrm{Im}\Omega.(M)$, that is, $N\simeq {}_fM$ as bicomodules for
some $f\in \aut{C}$.
\end{corollary}

\section{The Aut-Pic Property}

Following \cite{Bolla:1984} and
\cite{Cuadra/Garcia/Torrecillas:2000}, we give the following

\begin{definition}
We say that a coring $\coring{C}$ has the \emph{right Aut-Pic
property} if the group morphism $\Omega$ of Theorem
\ref{coringexseq} is surjective. In such a case,
$\rout{\coring{C}}\simeq \rpic{\coring{C}}$.
\end{definition}

Bolla proved (\cite[p. 264]{Bolla:1984}) that every ring such that
all right progenerators (=finitely generated projective generators)
are free, has Aut-Pic. In particular, local rings (by using a
theorem of Kaplansky), principal right ideal domains (by
\cite[Corollary 2.27]{Lam:1999}), and polynomial rings
$k[X_1,\dots,X_n]$, where $k$ is a field (by Quillen-Suslin
Theorem), have Aut-Pic. He also proved that every basic
(semiperfect) ring has Aut-Pic (see \cite[Proposition
3.8]{Bolla:1984}). Moreover, it is well-known that a semiperfect
ring $A$ is Morita equivalent to its basic ring $eAe$ (see
\cite[Proposition 27.14]{Anderson/Fuller:1992}), and by Theorem
\ref{coringexseq}, $\pic{A}\simeq \out{eAe}$. In
\cite{Cuadra/Garcia/Torrecillas:2000}, the authors gave several
interesting examples of coalgebras over fields having Aut-Pic. For
instance, they proved that every basic coalgebra has Aut-Pic. On the
other hand, we know \cite[Corollary 2.2]{Chin/Montgomery:1995} that
given a coalgebra over a field $C$ is Morita-Takeuchi equivalent to
a basic coalgebra $C_0$, and by Theorem \ref{coringexseq},
$\pic{C}\simeq \out{C_0}$.

\medskip
Of course all of the examples of rings and coalgebras over fields
having Aut-Pic mentioned in \cite{Bolla:1984} and
\cite{Cuadra/Garcia/Torrecillas:2000} are examples of corings having
Aut-Pic. In order to give others examples of corings satisfying the
Aut-Pic property we need the next lemma which is a generalization of
\cite[Proposition 4.1]{Nastasescu/Torrecillas/VanOystaeyen:2000}.

\medskip
We recall from \cite{Nastasescu/Torrecillas/VanOystaeyen:2000}, that
an object $M$ in an additive category $\cat{C}$ has the
\emph{invariant basis number property} (IBN for short) if $M^n\simeq
M^m$ implies $n=m$. For example every non-zero finitely generated
projective module over a semiperfect ring has IBN (from
\cite[Theorem 27.11]{Anderson/Fuller:1992}). A bicomodule $M \in
{}_B\rcomod{\coring{C}}$ is said to be
$(B,\coring{C})$-\emph{quasi-finite} \cite{Brzezinski/Wisbauer:2003}
if the functor $- \tensor{B} M : \rmod{B} \rightarrow
\rcomod{\coring{C}}$ has a left adjoint $\cohom{\coring{C}}{M}{-} :
\rcomod{\coring{C}} \rightarrow \rmod{B}$, and we call it the
\emph{cohom functor}. $M$ is said to be
$(B,\coring{C})$-\emph{injector} if the functor $- \tensor{B} M :
\rmod{B} \rightarrow \rcomod{\coring{C}}$ preserves injective
objects. Of course, if $B$ is a QF ring and $M$ is
$(B,\coring{C})$-injector then $M$ is injective in
$\rcomod{\coring{C}}$.

\begin{lemma}\label{IBN}
Let $\coring{C}$ be an $A$-coring such that $_A\coring{C}$ is flat.
If $0\neq M\in{}_B\rcomod{\coring{C}}$ is quasi-finite, such that
\begin{enumerate}[(a)]
\item $B$ is semisimple and $\rcomod{\coring{C}}$ is locally finitely
generated; or
\item $B$ is semiperfect and there is a finitely generated
projective $M_0$ in $\rcomod{\coring{C}}$ such that
$\hom{\coring{C}}{M_{0}}{M}\neq 0$. (the last condition is fulfilled
if $\coring{C}$ is cosemisimple (= $\rcomod{\coring{C}}$ is a
discrete spectral category, see
\cite{Stenstrom:1975,Brzezinski/Wisbauer:2003}, or if $A$ is right
artinian, $_A\coring{C}$ is projective, and $\coring{C}$ is
semiperfect, see \cite[Theorem 3.1]{Caenepeel/Iovanov:2005}.);
\end{enumerate}
then $M$ has IBN as a $B-\coring{C}$-bicomodule.

In particular, if $A$ is semisimple, or $A$ is semiperfect and there
is a non-zero finitely generated projective in
$\rcomod{\coring{C}}$, then $0\neq \coring{C}$ has IBN as a
$A-\coring{C}$-bicomodule.
\end{lemma}

\begin{proof}
We prove at the same time the two statements. Let $M_0$ be a
finitely generated subcomodule of $M$ (resp. a finitely generated
projective in $\rcomod{\coring{C}}$ such that
$\hom{\coring{C}}{M_{0}}{M}\neq 0$). We have
$\cohom{\coring{C}}{M}{M_{0}}^*=\hom{B}{\cohom{\coring{C}}{M}{M_{0}}}{B}\simeq
\hom{\coring{C}}{M_{0}}{M}$ as left $B$-modules. Since the functor
$-\tensor{B}{M}:\rmod{B}\to\rcomod{\coring{C}}$ is exact and
preserves coproducts, the cohom functor
$\cohom{\coring{C}}{\Lambda}{-}:\rcomod{\coring{C}}\to\rmod{B}$
preserves finitely generated (resp. finitely generated projective)
objects. In particular, $\cohom{\coring{C}}{M}{M_{0}}$ is a finitely
generated (resp. finitely generated projective) right $B$-module.
Hence the left $B$-module $\hom{\coring{C}}{M_{0}}{M}$ is so.

Now suppose that $M^m\simeq M^n$. Since the functor
$\hom{\coring{C}}{M_{0}}{-}:{}_B\rcomod{\coring{C}}\to {}\lmod{B}$
is $k$-linear, then $\hom{\coring{C}}{M_{0}}{M^m}\simeq
\hom{\coring{C}}{M_{0}}{M^n}$ as left $B$-bimodules. Hence
$\hom{\coring{C}}{M_{0}}{M}^m \\ \simeq
\hom{\coring{C}}{M_{0}}{M}^n$ as left $B$-bimodules. Finally, by
\cite[Theorem 27.11]{Anderson/Fuller:1992}, $m=n$.

For the particular case we take $M=\coring{C}$.
\end{proof}

\begin{corollary}
Let $C$ be a $k$-coalgebra such that $_kC$ is flat. If $0\neq
M\in{}\rcomod{C}$ is quasi-finite, such that
\begin{enumerate}[(a)]
\item $k$ is semisimple; or
\item $k$ is semiperfect and there is a finitely generated
projective $M_0$ in $\rcomod{C}$ such that $\hom{C}{M_{0}}{M}\neq
0$;
\end{enumerate}
then $M$ has IBN as a ${C}$-comodule.

In particular, if $k$ is semisimple, or $k$ is semiperfect and there
is a non-zero finitely generated projective in $\rcomod{C}$, then
$0\neq C$ has IBN as an $C$-bicomodule.
\end{corollary}

\begin{proposition}
Let $\coring{C}\neq0$ be an $A$-coring such that $A$ is semisimple,
or $A$ is semiperfect and there is a non-zero finitely generated
projective in $\rcomod{\coring{C}}$. If every
$\coring{C}$-bicomodule which is $(A,\coring{C})$-injector, is
isomorphic to $\coring{C}^{(I)}$ as $A-\coring{C}$-bicomodules for
some set $I$, then the coring $\coring{C}$ has right Aut-Pic.
\end{proposition}

\begin{proof}
The proof is analogous to that of \cite[Proposition
2.4]{Cuadra/Garcia/Torrecillas:2000} and the proof of Bolla
(\cite[p. 264]{Bolla:1984}) of the fact that every ring such that
all left progenerators are free has Aut-Pic.

Let $M$ be a right invertible $\coring{C}$-bicomodule. By
assumptions, $M\simeq \coring{C}^{(I)}$ as
$A-\coring{C}$-bicomodules for some set $I$. Let $M_0$ be a finitely
generated (resp. finitely generated projective) right comodule.
Since $M_0$ is a small object (see \cite{Popescu:1973}), then the
$k$-linear functor
$\hom{\coring{C}}{M_{0}}{-}:\rcomod{\coring{C}}\to {}\rmod{k}$
preserves coproducts. Therefore, $\hom{\coring{C}}{M_{0}}{M}\simeq
\hom{\coring{C}}{M_{0}}{\coring{C}}^{(I)}$ as left $A$-modules.
Hence I is a finite set, and $M\simeq \coring{C}^{(n)}$ as
$A-\coring{C}$-bicomodules for some $n\geq1$. Let $N$ be a
bicomodule such that $(N)$ is the inverse of $(M)$ in
$\rpic{\coring{C}}$, then
$$\coring{C}\simeq M\cotensor{\coring{C}}N\simeq
\coring{C}^{(n)}\cotensor{\coring{C}}N\simeq N^{(n)}$$ as
$A-\coring{C}$-bicomodules. On the other hand, $N\simeq
\coring{C}^{(m)}$ as $A-\coring{C}$-bicomodules for some set
$m\geq1$. It follows that $\coring{C}\simeq \coring{C}^{(nm)}$ as
$A-\coring{C}$-bicomodules. By Lemma \ref{IBN}, $nm=1$ and then
$M\simeq \coring{C}$ as $A-\coring{C}$-bicomodules. Finally, from
Proposition \ref{oneside}, $M\simeq {}_f\coring{C}$ as bicomodules
for some $f=(\varphi,1_A)\in \aut{\coring{C}}$.
\end{proof}

\begin{corollary}
Let $C$ be a $k$-coalgebra such that $k$ is semisimple, or $k$ is a
QF ring and there is a non-zero finitely generated projective in
$\rcomod{C}$. If every right injective $C$-comodule is free, then
the coalgebra $C$ has right Aut-Pic.
\end{corollary}

The following result allow us to simplify the computation of the
Picard group of some interesting corings.

\begin{proposition}\label{examples}
\begin{enumerate}[(1)] \item Let $_B\Sigma_A$ be a bimodule
such that $\Sigma_A$ is finitely generated and projective and
$_B\Sigma$ is faithfully flat. Then
$\rpic{\Sigma^*\tensor{B}\Sigma}\simeq \pic{B}$, where
$\Sigma^*\tensor{B}\Sigma$ is the comatrix coring
\cite{ElKaoutit/Gomez:2003} associated to $\Sigma$. If moreover $B$
has Aut-Pic then $\rpic{\Sigma^*\tensor{B}\Sigma}\simeq \out{B}$.
\\ In particular, For a $k$-algebra $A$ and $n\in \mathbb{N}$, we have
$\rpic{M^c_n(A)}\simeq \pic{A}$, where $M^c_n(A)$ is the
$(n,n)$-matrix coring over $A$ defined in
\cite[17.7]{Brzezinski/Wisbauer:2003}. If moreover $A$ has Aut-Pic
then $\rpic{M^c_n(A)}\simeq \out{A}$.
\item Let $\coring{C}$ be an $A$-coring such that $\coring{C}_A$ is
flat, and $R$ the opposite algebra of $^*\coring{C}$. If
$_A\coring{C}$ is finitely generated projective (e.g. if
$\coring{C}$ is a Frobenius ring, see
\cite{Brzezinski/Wisbauer:2003}), then $\rpic{\coring{C}}\simeq
\pic{R}$.
 \item If $\coring{C}$ is an $A$-coring such that the category
 $\rcomod{\coring{C}}$ has a finitely generated projective
 generator $U$, then $\rpic{\coring{C}}\simeq
\pic{\operatorname{End}_\coring{C}(U)}$.
\end{enumerate}
\end{proposition}

\begin{proof}
$(1)$ From \cite[Theorem 3.10]{ElKaoutit/Gomez:2003}, we have an
equivalence of categories $-\tensor{B}\Sigma:{}\rmod{B}\to
\rcomod{\Sigma^*\tensor{B}\Sigma}$. Theorem \ref{coringexseq}
achieves then the proof. Now we prove the particular case. Let $A$
be a $k$-algebra and $n\in \mathbb{N}$. If we take $\Sigma=A^n$, the
comatrix coring $\Sigma^*\tensor{B}\Sigma$ can be identified with
$M^c_n(A)$.

$(2)$ From \cite[Lemma 4.3]{Brzezinski:2002}, the categories
$\rcomod{\coring{C}}$ and $\rmod{R}$ are isomorphic to each other.
It is enough to apply Theorem \ref{coringexseq}.

$(3)$ This follows immediately from an alternative of
Gabriel-Popescu's Theorem (see \cite[p. 223]{Stenstrom:1975}) and
Theorem \ref{coringexseq}.
\end{proof}

\section{Application to the Picard group of $gr-(A,X,G)$}

In this Section we adopt the notations of \cite{DelRio:1992} and
\cite{Caenepeel/Militaru/Zhu:2002}.

We recall from \cite{Brzezinski/Wisbauer:2003}, that a
\emph{right-right entwining structure over} $k$ is a triple
$(A,C,\psi)$, where $A$ is a $k$-algebra, $C$ is a $k$-coalgebra,
and $\psi:C\otimes A\to A\otimes C$ is a $k$-linear map, such that
\begin{enumerate}[(ES1)]
\item $\psi\circ(1_C\otimes m)=(m\otimes 1_C)\circ(1_A\otimes
\psi)\circ(\psi\otimes 1_A)$, or equivalently, for all $a,b\in A$,
$c\in C$, $\sum(ab)_\psi\otimes c^\psi=\sum a_\psi b_\Psi\otimes
c^{\psi\Psi}$,
\item $(1_a\otimes \Delta)\circ \psi=(\psi\otimes1_C)\circ (1_C\otimes
\psi)\circ(\Delta\otimes 1_A)$, or equivalently, for all $a\in A$,
$c\in C$, $\sum a_\psi\otimes \Delta(c^\psi)=\sum
a_{\psi\Psi}\otimes c_{(1)}^\Psi\otimes c_{(2)}^\psi$,
\item $\psi\circ(1_C\otimes\eta)=\eta\otimes1_C$, or equivalently, for
all $c\in C$, $\sum(1_A)_\psi\otimes c^\psi=1_A\otimes c$,
\item $(1_A\otimes\epsilon)\circ\psi=\epsilon\otimes1_A$, or
equivalently, for all $a\in A$, $c\in C$, $\sum
a_\psi\epsilon(c^\psi)=a\epsilon(c)$.
\end{enumerate}
where $m$ and $\eta$ are respectively the multiplication and the
unit maps of $A$, and $\psi(c\otimes a)=a_\psi\otimes
c^\psi=a_\Psi\otimes c^\Psi.$

A \emph{morphism} $(A,C,\psi)\to(A',C',\psi')$ is a pair
$(\alpha,\gamma)$ with  $\alpha:A\to A'$ is a morphism of algebras,
and $\gamma:C\to C'$ is a morphism of coalgebras such that
$(\alpha\otimes\gamma)\circ\psi=\psi'\circ(\gamma\otimes\alpha)$ or
equivalently, for all $a\in A$, $c\in C$,
$\sum\alpha(a_\psi)\otimes\gamma(c^\psi)=\sum\alpha(a)_{\psi'}\otimes\gamma(c)^{\psi'}$.
We denote this category by $\mathbb{E}_\bullet^\bullet(k)$. Notice
that $(\alpha,\gamma)\in \mathbb{E}_\bullet^\bullet(k)$ is an
isomorphism if and only if $\alpha$ and $\gamma$ are bijective.

Again we recall from \cite{Brzezinski/Wisbauer:2003}, that a
\emph{right-right entwined modules} over a right-right entwining
structure $(A,C,\psi)$ is a $k$-module $M$ which is a right
$A$-module with multiplication $\psi_M$, and a right $C$-comodule
with comultiplication $\rho_M$ such that
$(\psi_M\otimes1_C)\circ(1_M\otimes\psi)\circ(\rho_M\otimes1_A)=\rho_M\circ\psi_M$,
or equivalently, for all $m\in M$ and $a\in A$, $\rho_M(ma)=\sum
m_{(0)}a_\psi\otimes m_{(1)}^\psi$. A \emph{morphism} between
entwined modules is a morphism of right $A$-modules and right
$C$-comodules at the same time. We denote this category by
$\mathcal{M}(\psi)_A^C$.

The following result due to Takeuchi shows that entwining structures
and entwined modules are very related to corings and comodules over
corings respectively. For a better understanding, we add a sketch of
the proof.

\begin{theorem}
\begin{enumerate}[(a)]
\item Let $A$ be an algebra, $C$ be a $k$-coalgebra, and let
$\psi:C\otimes A\to A\otimes C$ be a $k$-linear map. Obviously
$A\otimes C$ has a structure of left $A$-module by $b(a\otimes
c)=ba\otimes c$ for $a,b\in A$ and $c\in C$. Define the right
$A$-module action on $A\otimes C$, $\xymatrix{\psi_{A\otimes
C}^r:A\otimes C\otimes A \ar[r]^-{A\otimes\psi} & A\otimes A\otimes
C \ar[r]^-{m\otimes C} & A\otimes C}$, that is, $(a\otimes
c)b=a\psi(c\otimes b)$ for $a,b\in A$ and $c\in C$. Define also
$$\xymatrix{\Delta:A\otimes C\ar[r]^-{A\otimes \Delta_C} & A\otimes
C\otimes C\simeq (A\otimes C)\tensor{A}(A\otimes C)}$$ (for every
$a\in A,c\in C,$ $\Delta(a\otimes c)=\sum(a\otimes
c_{(1)})\tensor{A}(1_A\otimes c_{(2)})$, where $\Delta_C(c)=\sum
c_{(1)}\otimes c_{(2)}$, and
$$\xymatrix{\epsilon:A\otimes C\ar[r]^-{A\otimes \epsilon_C} &
A\otimes k\simeq A}$$ ($\epsilon(a\otimes c)=a\epsilon_C(c)$). Thus,
$(A\otimes C,\Delta,\epsilon)$ is an $A$-coring if and only if
$(A,C,\psi)$ is an entwining structure.
\item Let $(A,C,\psi)$ be an entwining
structure, $M$ be a right $A$-module with the action
$\psi_M:M\otimes A\to M$, and let $\rho_M:M\to M\otimes C$ be a
$k$-linear map. Define $$\xymatrix{\rho'_M:M\ar[r]^-{\rho_M} &
M\otimes C\simeq M\tensor{A}(A\otimes C)}$$ (for every $m\in M,$
$\rho'_M(m)=\sum m_{(0)}\tensor{A}(1_A\otimes m_{(1)})$, where
$\rho_M(m)=\sum m_{(0)}\otimes m_{(1)}$). Then, $(M,\rho'_M)$ is a
right $A\otimes C$-comodule if and only if $(M,\psi_M,\rho_M)$ is a
right-right entwined module over $(A,C,\psi)$. Hence, the category
of right $A\otimes C$-comodules, $\rcomod{A\otimes C}$ is isomorphic
to the category of right-right entwined modules over $(A,C,\psi)$,
$\mathcal{M}(\psi)_A^C$.
\end{enumerate}
\end{theorem}

\begin{proof}
$(a)$ It is easy to verify that\begin{itemize}
\item For all $a,b\in A,c\in C$, $(1\otimes
c)(ab)=((1\otimes c)a)b \Longleftrightarrow (\textrm{ES}1)$.
\item For all $c\in
C$, $(1\otimes c)1=1\otimes c \Longleftrightarrow (\textrm{ES}3)$.
\item For all $a\in A,c\in C$, $\Delta((1\otimes c)a)=\Delta(1\otimes
c)a \Longleftrightarrow (\textrm{ES}2)$. \item For all $a\in A,c\in
C$, $\epsilon((1\otimes c)a)=\epsilon(1\otimes c)a
\Longleftrightarrow (\textrm{ES}4)$.
\end{itemize}
Moreover, the coassociativity of $\Delta$ and the counit property of
$\epsilon$ follow from that of $\Delta_C$ and $\epsilon_C$. Hence
$(a)$ follows.

$(b)$ It is easy to verify that\begin{itemize}
\item $\rho'_M$ is $A$-linear if and only if for all $m\in M$ and $a\in
A$, $\rho_M(ma)=\sum m_{(0)}a_\psi\otimes m_{(1)}^\psi$.
\item $\rho'_M$ is coassociative if and only if $\rho_M$ is
coassociative.
\item The counit property of $\epsilon$ holds if and only if that of
$\epsilon_C$ holds.
 Hence $(b)$ follows.
\end{itemize}
\end{proof}
 It is easy to verify that if
$(\alpha,\gamma):(A,C,\psi)\rightarrow(A',C',\psi')$ is a morphism
of entwined structures, then $(\alpha\otimes\gamma,\alpha):A\otimes
C\to A\otimes C$ is a morphism of corings. Hence we have a functor
$$F:\mathbb{E}_\bullet^\bullet(k)\to\mathbf{Crg}_k.$$
By Theorem \ref{coringexseq}, we obtain

\begin{proposition}\label{entwining}
For an entwining structure $(A,C,\psi)$ there is an exact sequence
$$\xymatrix{1\ar[r] & \operatorname{Ker}(\Omega\circ F)\ar[r]&
\aut{(A,C,\psi)} \ar[r]^-{\Omega\circ F} & \rpic{A\otimes C}}\simeq
\rpic{\mathcal{M}(\psi)_A^C},$$ and $\operatorname{Ker}(\Omega\circ
F)$ is the set of $(\alpha,\gamma)\in \aut{(A,C,\psi)}$ such that
there is $p\in (A\otimes C)^*$ satisfying $\sum
\big(\alpha(a)\otimes\gamma(c_{(1)})\big)p(1_A\otimes c_{(2)})=\sum
p(a\otimes c_{(1)})\otimes c_{(2)}$ for all $a\in A,c\in C$.
\end{proposition}

A right-right \emph{Doi-Koppinen structure or simply DK structure}
over $k$ \cite{Caenepeel/Militaru/Zhu:2002} is a triple $(H,A,C)$,
where $H$ is a bialgebra, $A$ is a right $H$-comodule algebra, and
$C$ is a right $H$-module coalgebra. A \emph{morphism of DK
structures} is a triple
$(\hbar,\alpha,\gamma):(H,A,C)\to(H',A',C')$, where $\hbar:H\to H'$,
$\alpha:A\to A'$, and $\gamma:C\to C'$ are respectively a bialgebra
morphism, an algebra morphism, and a coalgebra morphism such that
$\rho_{A'}(\alpha(a))=\alpha(a_{(0)})\otimes \hbar(a_{(1)})$ and
$\gamma(ch)=\gamma(c)\hbar(h)$, for all $a\in A, c\in C, h\in H$.
This yields a category which we denote by
$\mathbb{DK}_\bullet^\bullet(k)$. Moreover, we have a faithful
functor $G:\mathbb{DK}_\bullet^\bullet(k)\to
\mathbb{E}_\bullet^\bullet(k)$ defined by $G((H,A,C))=(A,C,\psi)$
with $\psi:C\otimes A\to A\otimes C$ and $\psi(c\otimes
a)=a_{(0)}\otimes ca_{(1)}$, and
$G((\hbar,\alpha,\gamma))=(\alpha,\gamma)$. (see \cite[Proposition
17]{Caenepeel/Militaru/Zhu:2002}.) Notice that
$(\hbar,\alpha,\gamma)\in \mathbb{DK}_\bullet^\bullet(k)$ is an
isomorphism if and only if $\hbar$, $\alpha$ and $\gamma$ are
bijective. The category of right \emph{Doi-Koppinen-Hopf modules}
over the right-right DK structure $(H,A,C)$ is exactly the category
of $\mathcal{M}(\psi)_A^C$, and it is denoted by
$\mathcal{M}(H)_A^C$.

\medskip
Hence Proposition \ref{entwining} yields the following

\begin{proposition}\label{DK}
For a DK structure $(H,A,C)$ there is an exact sequence
$$\xymatrix{1\ar[r] & \operatorname{Ker}(\Omega\circ F\circ G)\ar[r]&
\aut{(H,A,C)}\ar[r]^-{\Omega\circ F\circ G} & \rpic{A\otimes
C}\simeq\rpic{\mathcal{M}(H)_A^C}},$$ and
$\operatorname{Ker}(\Omega\circ F\circ G)$ is the set of
$(\hbar,\alpha,\gamma)\in \aut{(H,A,C)}$ such that there is $p\in
(A\otimes C)^*$ satisfying $\sum
\big(\alpha(a)\otimes\gamma(c_{(1)})\big)p(1_A\otimes c_{(2)})=\sum
p(a\otimes c_{(1)})\otimes c_{(2)}$ for all $a\in A,c\in C$.
\end{proposition}

Finally, we consider the category of right modules graded by a
$G$-set, $gr-(A,X,G)$, where $G$ is a group and $X$ is a right
$G$-set. This category is introduced and studied in
\cite{Nastasescu/Raianu/VanOystaeyen:1990}. A study of the graded
ring theory can be found in the recent book
\cite{Nastasescu/VanOystaeyen:2004}. This category is equivalent to
a category of modules over a ring if $X=G$ and $A$ is strongly
graded (by Dade Theorem \cite[Theorem
3.1.1]{Nastasescu/VanOystaeyen:2004}), or if $X$ is a finite set (by
\cite[Theorem 2.13]{Nastasescu/Raianu/VanOystaeyen:1990}). But there
is an example of a category of graded modules which is not
equivalent to a category of modules, see \cite[Remark
2.4]{Menini/Nastasescu:1988}.

Let $G$ be a group, $X$ a right $G$-set, and $A$ be a $G$-graded
$k$-algebra. We know that $(kG,A,C)$ is a DK structure, with $kG$ is
a Hopf algebra. Then we have, $(A,kX,\psi)$ is an entwining
structure where $\psi:kX\otimes A\rightarrow A\otimes kX$ is the map
defined by $\psi(x\otimes a_{g})=a_{g}\otimes xg$ for all $x\in
X,g\in G,a_g\in A_g$. From the above considerations, we have an
$A$-coring, $A\otimes kX$. The comultiplication and the counit maps
of the coring $A\otimes kX$ are defined by:
$$\Delta_{A\otimes kX}(a\otimes x)=(a\otimes
x)\tensor{A}(1_A\otimes x), \quad \epsilon_{A\otimes kX}(a\otimes
x)=a \quad (a\in A,x\in X).$$ We know that this coring is
coseparable (see \cite[Section 6]{Zarouali:2004}). From \cite[\S
4.6]{Caenepeel/Militaru/Zhu:2002}, $gr-(A,X,G)\simeq
\mathcal{M}(kG)_A^{kX}$. Moreover, the left $A$-module $A\otimes kX$
is free with the basis $\{1_A\otimes x\mid x\in X\}$. It is easy to
see that the last family is a basis of the right $A$-module
$A\otimes kX$ (since $a_g\otimes x=(1_A\otimes xg^{-1})a_g$ for
every $g\in G,a_g\in A_g,x\in X$, and since $\psi$ is an
isomorphism.) Each $p\in (A\otimes kX)^*$ is entirely determined by
the data of $p(1_A\otimes x)$, for all $x\in X$. The same thing
holds for $\varphi$, where $(\varphi,\rho)\in \aut{A\otimes kX}$.

\medskip
By Theorem \ref{coringexseq} we get

\begin{proposition}
\begin{enumerate}[(a)]\item We have, $p\in (A\otimes kX)^*$ is
invertible
 if and only if there exists $q\in (A\otimes kX)^*$ such that
 $\sum_{h\in G}q(1_A\otimes xh^{-1})p(1_A\otimes
x)_h=1_A=\sum_{h\in G}p(1_A\otimes xh^{-1})q(1_A\otimes x)_h$ for
every $x\in X$.
\item We have an exact sequence
$$\xymatrix{1\ar[r] & \operatorname{Ker}(\Omega)\ar[r]&
\aut{A\otimes kX}\ar[r]^-\Omega & \pic{A\otimes kX}\simeq
\pic{gr-(A,X,G)}},$$ and $(\varphi,\rho)\in
\operatorname{Ker}(\Omega)$ if and only if there exists $p\in
(A\otimes kX)^*$ invertible such that
\begin{enumerate}[(i)] \item for
all $x\in X$, $a_y^xp(1_A\otimes x)_h=0$ for all $y\in X,h\in G$
such that $yh\neq x$, where $\varphi(1_A\otimes x)=\sum_{y\in
X}a_y^x\otimes y$, and \item for all $x\in X,a\in A$, $p(a\otimes
x)=\rho(a)p(1_A\otimes x).$
\end{enumerate}
\end{enumerate}
\end{proposition}

\begin{proof}
We have, $p\in (A\otimes kX)^*$ is invertible if and only if there
exists $q\in (A\otimes kX)^*$ such that
\begin{eqnarray}
q\big(p(a\otimes x)(1_A\otimes x)\big)=\epsilon(a\otimes x)=a \\
p\big(q(a\otimes x)(1_A\otimes x)\big)=\epsilon(a\otimes x)=a,
\end{eqnarray}
for all $a\in A, x\in X$. (from
\cite[17.8]{Brzezinski/Wisbauer:2003}.) On the other hand,
\begin{eqnarray*}
q\big(p(a_g\otimes x)(1_A\otimes x)\big) &=& q\big(p(1_A\otimes
xg^{-1})(a_g\otimes x)\big) \\
&=& q\big(\sum_{h\in G}(b_ha_g)\otimes
x\big),\quad\textrm{where}\quad p(1_A\otimes xg^{-1})=\sum_{h\in
G}b_h \\
&=& \sum_{h\in G}q\big(1_A\otimes x(hg)^{-1}\big)(b_ha_g), \quad
\textrm{since}\quad b_ha_g\in A_{hg} \\
&=& \sum_{h\in G}q(1_A\otimes xg^{-1}h^{-1})p(1_A\otimes
xg^{-1})_ha_g,
\end{eqnarray*}
for all $x\in X,g\in G, a_g\in A_g$. Hence
\begin{eqnarray*}
(2)&\Longleftrightarrow& q\big(p(a_g\otimes x)(1_A\otimes x)\big)=
a_g\quad \textrm{for all}\, x\in X,g\in G, a_g\in A_g \\
&\Longleftrightarrow& \sum_{h\in G}q(1_A\otimes xh^{-1})p(1_A\otimes
x)_h=1_A \quad \textrm{for all}\, x\in X.
\end{eqnarray*}
By symmetry, $$(3)\Longleftrightarrow \sum_{h\in G}p(1_A\otimes
xh^{-1})q(1_A\otimes x)_h=1_A \quad \textrm{for all}\, x\in X.$$

By Theorem \ref{coringexseq}, $(\varphi,\rho)\in
\operatorname{Ker}(\Omega)$ if and only if there exists $p\in
(A\otimes kX)^*$ invertible with $\varphi(a\otimes x)p(1_A\otimes
x)=p(a\otimes x)\otimes x$ for all $a\in A,x\in X$. On the other
hand,
\begin{eqnarray*}
\varphi(a\otimes x)p(1_A\otimes x)&=& \rho(a)\varphi(1_A\otimes
x)p(1_A\otimes x) \\
&=& \rho(a)\sum_{y\in X}\sum_{h\in G}(a_y^x\otimes
y)b_h,\,\textrm{where}\, \varphi(1_A\otimes x)=\sum_{y\in
X}a_y^x\otimes y,\, p(1_A\otimes x)=\sum_{h\in G}b_h \\
&=& \rho(a)\sum_{y\in X}\sum_{h\in G}a_y^xb_h\otimes yh.
\end{eqnarray*}
Since $\sum_{y\in X}a_y^x=\epsilon\circ\varphi(1_A\otimes
x)=\rho\circ\epsilon(1_A\otimes x)=1_A$, $(b)$ follows.
\end{proof}

Now, let $f:G\rightarrow G'$ be a morphism of groups, $X$ a right
$G$-set, $X'$ a right $G'$-set, and $\varphi:X\rightarrow X'$ a map.
Let $A$ be a $G$-graded $k$-algebra, $A'$ a $G'$-graded $k$-algebra,
and $\alpha:A\rightarrow A'$ a morphism of algebras. We have,
$\hbar:kG\rightarrow kG'$ defined by $\hbar(g) =f(g)$ for each $g\in
G$, is a morphism of Hopf algebras, and $\gamma:kX\rightarrow kX'$
defined by $\gamma(x) =\varphi(x)$ for each $x\in X$, is a morphism
of coalgebras. It is easy to show that $(\hbar,\alpha,\gamma)$ is a
morphism of DK structures if and only if
\begin{eqnarray}
\varphi(xg) =\varphi(x)f(g)\quad  \textrm{for all}\, g\in G,x\in X \\
\alpha( A_{g}) \subset A_{f(g)}' \quad\textrm{for all}\, g\in G.
\end{eqnarray}

We define a category as follows. The objects are the triples
$(G,X,A)$, where $G$ is a group, $X$ is a right $G$-set, and $A$ is
a graded $k$-algebra. The morphisms are the triples
$(f,\varphi,\alpha)$, where $f:G\rightarrow G'$ is a morphism of
groups, $\varphi:X\rightarrow X'$ is a map, and $\alpha:A\rightarrow
A'$ is a morphism of algebras such that $(6),(7)$ hold. We denote
this category by $\mathbb{G}^r(k)$. There is a faithful functor
$$H:\mathbb{G}^r(k)\to \mathbb{DK}_\bullet^\bullet(k)$$ defined by
$H((G,X,A))=(kG,A,kX)$ and
$H((f,\varphi,\alpha))=(\hbar,\alpha,\gamma)$. We obtain then the
next result. We think that our subgroup of $\pic{gr-(A,X,G)}$ is
more simple than the Beattie-Del R\'io subgroup (see \cite[\S
2]{Beattie/Delrio:1996}).

\begin{proposition}\label{graded}
There is an exact sequence
$$\xymatrix{1\ar[r] & \operatorname{Ker}(\Omega\circ F\circ G\circ
H)\ar[r]& \aut{(G,X,A)}\ar[rr]^-{\Omega\circ F\circ G\circ H} & &
\pic{A\otimes kX}\simeq \pic{gr-(A,X,G)}},$$ and
$\operatorname{Ker}(\Omega\circ F\circ G\circ H)$ is the set of
$(f,\varphi,\alpha)\in \aut{(G,X,A)}$ such that there is $p\in
(A\otimes kX)^*$ satisfying $p(a\otimes x)=\alpha(a)p(1_A\otimes x)$
for all $x\in X,a\in A$; and $p(1_A\otimes x)_h=0$ for all $x\in X,
h\in G$ such that $\varphi(x)h\neq x$.
\end{proposition}

\section*{Acknowledgements}
I would like to thank my advisor, Professor Jos\'e
G\'omez-Torrecillas, for the subject of this paper. I am also very
grateful to Professor Blas Torrecillas for the helpful discussion.

\noindent Department of Algebra, University of Granada,\\
E18071 Granada, Spain \\
E-mail address: \textsf{zaroual@correo.ugr.es}

\end{document}